\documentclass[]{nyjm}

\pdfoutput=1

\usepackage{amssymb}
\usepackage{graphicx}
\usepackage{amscd}
\usepackage{hyperref}
\usepackage[usenames,dvipsnames]{color}
\hypersetup{
  colorlinks,
  citecolor=Blue,
  linkcolor=Blue,
  urlcolor=Blue}

\title{CO-HIGGS BUNDLES ON $\PP^1$}
\author{STEVEN RAYAN}
\address{Dept. of Mathematics, Univ. of Toronto, 40 St. George St., Toronto, ON, Canada, M5S 2E4.}
\email{\href{mailto:rayan@math.toronto.edu}{rayan@math.toronto.edu}}
\urladdr{This paper is available via \href{http://nyjm.albany.edu/j/2013/19-42.html}{\ttfamily http://nyjm.albany.edu/j/2013/19-42.html}.}

\thanks{Parts of this work were funded by the Commonwealth Scholarship Plan and the Natural Sciences and Engineering Research Council of Canada.}



\newtheorem{theorem}{Theorem}[section]

\newtheorem{proposition}[theorem]{Proposition}

\theoremstyle{definition}
\newtheorem{definition}[theorem]{Definition}

\newtheorem{remark}[theorem]{Remark}
\newtheorem*{note}{Notation}

\newcommand{\loclab}{\label{\arabic{section}.\arabic{theorem}.\theenumi}}

\numberwithin{equation}{section}

\newcommand{\bdoc}{\begin{document}}
\newcommand{\edoc}{\end{document}}

\newcommand{\bcent}{\begin{center}}
\newcommand{\ecent}{\end{center}}

\newcommand{\benum}{\begin{enumerate}}
\newcommand{\eenum}{\end{enumerate}}
\newcommand{\bitem}{\begin{itemize}}
\newcommand{\eitem}{\end{itemize}}

\newcommand{\btab}{\begin{tabular}}
\newcommand{\etab}{\end{tabular}}
\newcommand{\beqn}{\begin{eqnarray}}
\newcommand{\eeqn}{\end{eqnarray}}

\newcommand{\bmath}{\begin{math}}
\newcommand{\emath}{\end{math}}

\newcommand{\noin}{\noindent}

\providecommand{\tb}[1]{\textbf{#1}}

\newcommand{\bsh}{\backslash}

\newcommand{\ds}{\displaystyle}

\newcommand{\ub}{\underbrace}
\newcommand{\ob}{\overbrace}

\providecommand{\F}[1]{\mathbb{#1}}  

\newcommand{\mb}{\mathbf}

\newcommand{\FF}{\F F}
\providecommand{\Fn}[1]{\FF_{#1}}
\newcommand{\Fp}{\Fn{p}}
\newcommand{\Fq}{\Fn{q}}
\newcommand{\Fpm}{\Fn{p^m}}
\newcommand{\Fpn}{\Fn{p^n}}
\newcommand{\Fpr}{\Fn{p^r}}

\newcommand{\ZZ}{\mathbb{Z}} 
\providecommand{\Zn}[1]{\ZZ_{#1}}
\providecommand{\ZnZ}[1]{\ZZ/#1\ZZ}
\newcommand{\Zp}{\Zn{p}}
\newcommand{\ZpZ}{\ZnZ{p}}

\newcommand{\NN}{\F{N}}
\newcommand{\QQ}{\F{Q}}
\newcommand{\RR}{\F{R}}
\newcommand{\CC}{\F{C}}
\newcommand{\QQbar}{\overline{\QQ}}
\newcommand{\Zero}{\mathbb{00}}

\newcommand{\EE}{\F E}
\newcommand{\II}{\F I}
\newcommand{\KK}{\F K}
\newcommand{\MM}{\F M}
\newcommand{\XX}{\F X}
\newcommand{\PP}{\F P}
\newcommand{\FA}{\F A}
\newcommand{\LL}{\F L}
\newcommand{\HH}{\mathbb H} 
\newcommand{\FS}{\F S}
\newcommand{\FSig}{\F\Sig}
\newcommand{\FDel}{\F\Del}

\providecommand{\E}[1]{\hat{\F{#1}}}
\newcommand{\EC}{\E{C}}

\newcommand{\ind}{\mbox{ind}}

\newcommand{\fx}{f(x)}
\newcommand{\gx}{g(x)}

\newcommand{\x}{^\star}
\newcommand{\xs}{^{~\star}}
\providecommand{\U}[1]{\left(#1\right)\x}

\newcommand{\iso}{\simeq}
\newcommand{\plus}{\oplus}
\newcommand{\Plus}{\bigoplus}
\newcommand{\tensor}{\otimes}
\newcommand{\Tensor}{\bigotimes}
\newcommand{\inject}{\hookrightarrow}
\newcommand{\linject}{\hookleftarrow}
\newcommand{\surject}{\twoheadrightarrow}

\renewcommand{\ker}{\mbox{ker}\;}
\newcommand{\Imf}{\mbox{im}~\;}
\newcommand{\img}{\mbox{im}\;}
\newcommand{\Hom}{\mbox{Hom}}
\newcommand{\Sym}{\mbox{Sym}}
\newcommand{\End}{\mbox{End}\,}
\newcommand{\Endz}{\mbox{End}_0}
\newcommand{\Id}{\mbox{Id}}
\newcommand{\rk}{\mbox{rk}\,}
\newcommand{\Pic}{\mbox{Pic}}
\newcommand{\Jac}{\mbox{Jac}}
\newcommand{\ch}{\mbox{ch}\,}
\newcommand{\td}{\mbox{td}\,}

\providecommand{\Gal}[1]{\mbox{Gal}(#1)}
\providecommand{\GAL}[2]{\mbox{Gal}(#1/#2)}
\providecommand{\Sub}[1]{\mbox{Sub}(#1)}
\providecommand{\Lat}[1]{\mbox{Lat}(#1)}

\providecommand{\lst}[1]{\ds\left[\,\,#1\,\,\right]}

\newcommand{\dup}{d_\wedge}
\newcommand{\drt}{d_>}

\newcommand{\MF}{\mathfrak}

\newcommand{\Div}{~\lvert~}
\newcommand{\lcm}{\mbox{lcm}}

\providecommand{\leg}[2]{\left(\frac{#1}{#2}\right)}
\providecommand{\jac}[2]{\leg{#1}{#2}}
\providecommand{\qdc}[2]{\left[\frac{#1}{#2}\right]}

\newcommand{\w}{\omega}
\newcommand{\W}{\Omega}
\providecommand{\Cal}[1]{\mathcal{#1}}
\newcommand{\CL}{\Cal L}
\newcommand{\CO}{\Cal O}
\renewcommand{\O}{\Cal O}
\renewcommand{\o}{\Cal O}
\newcommand{\co}{\Cal O}
\newcommand{\CE}{\Cal E}
\newcommand{\CF}{\Cal F}
\newcommand{\CQ}{\Cal Q}
\newcommand{\CK}{\Cal K}
\newcommand{\CM}{\Cal M}
\newcommand{\CDee}{\Cal D}
\newcommand{\CP}{\Cal P}
\newcommand{\CN}{\Cal N}
\newcommand{\CS}{\Cal S}
\newcommand{\ClC}{\Cal C}

\newcommand{\Si}{\Sigma}

\providecommand{\Ok}[1]{\CO_{#1}}
\newcommand{\OK}{\Ok{K}}

\renewcommand{\epsilon}{\varepsilon}
\newcommand{\ep}{\varepsilon}

\providecommand{\abs}[1]{\left|#1\right|}
\providecommand{\norm}[1]{\lVert#1\rVert}

\newcommand{\di}{\partial}
\providecommand{\ddy}[1]{\ds\frac{d}{d #1}}
\providecommand{\didiy}[1]{\ds\frac{\di}{\di #1}}
\newcommand{\ddx}{\ddy{x}}
\newcommand{\didix}{\didiy{x}}
\providecommand{\v}[1]{\vec{#1}}

\newcommand{\CJ}{\Cal{J}}
\newcommand{\CA}{\Cal{A}}

\newcommand{\ihat}{\hat{\infty}}

\providecommand{\set}[1]{\left\{#1\right\}}

\providecommand{\ip}[2]{\left( #1,#2\right)}
\providecommand{\IP}[2]{\left\langle #1,#2\right\rangle}
\providecommand{\bra}[1]{\left\langle\left. #1\right|\right.}
\providecommand{\ket}[1]{\left.\left| #1\right.\right\rangle}
\providecommand{\bkv}[4]{\left\langle\left.\tb{#1} #2\right|\tb{#3} #4\right\rangle}
\providecommand{\bk}[2]{\bkv{}{#1}{}{#2}}
\providecommand{\bkvop}[5]{\left\langle\tb{#1} #2\left| #5\right|\tb{#3} #4\right\rangle}
\providecommand{\bkop}[3]{\bkvop{}{#1}{}{#2}{#3}}
\providecommand{\comm}[2]{\left[#1,#2\right]}
\providecommand{\expn}[1]{\left\langle #1\right\rangle}
\newcommand{\Del}{\Delta}
\newcommand{\Nab}{\nabla}
\newcommand{\Sig}{\Sigma}
\newcommand{\oline}{\overline}

\newcommand{\p}{\rho}
\providecommand{\sprod}[2]{\left\langle #1,#2\right\rangle}

\newcommand{\Ehat}{\overline E}
\newcommand{\Phihat}{\overline \Phi}
\newcommand{\phihat}{\overline \phi}

\newcommand{\QED}{\begin{flushright}\rule{2.5mm}{2.5mm}\\\end{flushright}}

\newenvironment{summary}{\noin\tb{Summary:~}}{}
\newenvironment{details}{\scriptsize\noin$\dagger$}{\normalsize}

\providecommand{\summary}[1]{\begin{summary}#1\end{summary}}
\providecommand{\details}[1]{\begin{details}#1\end{details}}

\begin{document}

\begin{abstract}  Co-Higgs bundles are Higgs bundles in the sense of Simpson, but with Higgs fields that take values in the tangent bundle instead of the cotangent bundle.  Given a vector bundle on $\PP^1$, we find necessary and sufficient conditions on its Grothendieck splitting for it to admit a stable Higgs field.  We characterize the rank-2, odd-degree moduli space as a universal elliptic curve with a globally-defined equation.  For ranks $r=2,3,4$, we explicitly verify the conjectural Betti numbers emerging from the recent work of Chuang, Diaconescu, Pan, and Mozgovoy on the ADHM formula.  We state the result for $r=5$.
\end{abstract}

\maketitle
\tableofcontents

\section{Introduction}

Let $X$ be an algebraic variety with cotangent bundle $T^*$.  A \emph{Higgs bundle} on $X$, in the sense of Simpson \cite{TSCH:88}, is a vector bundle $E\rightarrow X$ together with a \emph{Higgs field} $\phi\in H^0(X;(\End E)\tensor T^*)$ for which 
\[
\phi\wedge\phi=0\in H^0(X;(\End E)\tensor\wedge^2T^*).
\]  
Higgs bundles have been studied intensely, and appear naturally in areas of mathematics as diverse as string theory and number theory --- see \cite{BGG:07} for an overview.  

An alternative kind of Higgs bundle arises when we replace $T^*$ with $T$ in the definition of the Higgs field.  We call these objects \emph{co-Higgs bundles}.  They are only beginning to attract interest; however, there are discussions related to them in \cite{NJH:10IIa,NJH:10Ia,SSR:11}.  One motivation for studying co-Higgs bundles comes from generalized geometry, because generalized holomorphic bundles on ordinary complex manifolds are precisely co-Higgs bundles \cite{MG:07}.

The purpose of this note is to characterize co-Higgs bundles over curves.  In this case, $\phi\wedge\phi=0$ is automatic.  From now on $X$ is a curve, by which we mean a nonsingular, connected, projective curve over $\CC$.
By vector bundle, we will always mean a holomorphic vector bundle.

We show that stability restricts our study to the projective line.  We then classify the vector bundles on $\PP^1$ admitting semistable Higgs fields by their splitting types, and use this classification to study explicitly
the odd-degree component of the rank-2 moduli space.  The main result is a global description of this smooth moduli space as the variety of solutions of an algebraic equation.  This equation is a universal one for the fibres of the associated Hitchin map, whose generic fibre in this case is a nonsingular elliptic curve.  
An immediate consequence of our description is that the Betti numbers of the moduli space are those of $S^2$.
In the even case, we characterize a section of the fibration by the splitting type of $E$. 

For $r=3$ and $r=4$ with odd degree, we use Morse theory to calculate the Betti numbers, verifying conjectural Betti numbers due to Chuang, Diaconescu, and Pan in \cite{CDP:11}, which was adapted to genus 0 by Mozgovoy in \cite{SM:12}.  We state the result for $r=5$ without proof, although the method of computation is described in $\S$\ref{SecBetti}.\enlargethispage{\baselineskip}

\begin{note} 
We denote the canonical line bundle of $X$ by $K$. Accordingly, the anticanonical line bundle  ---  equivalently the holomorphic tangent bundle  ---  is $K^*$.  As we agree that $\phi$ is always a $K^*$-valued endomorphism, there is no cause for confusion if we omit the parentheses around $\End E$ in $\phi\in H^0(X;(\End E)\tensor K^*)$.
\end{note}

\subsubsection*{Acknowledgements}  I thank Nigel Hitchin for pointing me to this topic and for his guidance.  I acknowledge Steven Bradlow, Jonathan Fisher, Peter Gothen, Marco Gualtieri, Tam\'as Hausel, Lisa Jeffrey, and Sergey Mozgovoy for enlightening discussions.  I thank Ruxandra Moraru for pointing out an error in a remark in the original manuscript, as well as the referee for suggesting corrections, clarifications, and a number of other improvements to the manuscript.

\section{Morphisms, stability, and $S$-equivalence}

The following notions carry over from Higgs bundles without modification.  A 
\emph{morphism} taking $(E,\phi)$ to $(E',\phi')$ is a commutative diagram\bcent\bmath\begin{CD}E @> \psi >> E'\\@V \phi VV @ VV \phi' V\\E\tensor K^* @> \psi\tensor1 >> E'\tensor K^*\end{CD}\emath\ecent\noin in which $\psi:E\rightarrow E'$ is a morphism of vector bundles.   The pairs $(E,\phi)$ and $(E',\phi')$ are 
\emph{isomorphic} or \emph{equivalent}
when we have such a diagram 
 in which 
$\psi$ is an isomorphism of bundles.  In particular, $(E,\phi)$ and $(E,\phi')$ are isomorphic if and only if there exists an automorphism $\psi$ of $E$ such that $\psi\phi\psi^{-1}=\phi'$.  

The appropriate stability condition for moduli of co-Higgs bundles on $X$ is Hitchin's slope-stability
condition, which he defined for Higgs bundles in \cite{NJH:86}.  Following his definition, we have:

\begin{definition}  A co-Higgs bundle $(E,\phi)$ over $X$ is \emph{(semi)stable} if
\beqn \frac{\deg U}{\mbox{rk}\,U} & < & \frac{\deg E}{\mbox{rk}\,E}\label{Stable}\eeqn\noin
(respectively, $\leq$) for each proper nonzero subbundle $U\subset E$ that is invariant under $\phi$ (meaning $\phi(U)\subseteq U\tensor K^*$).  The rational number 
\beqn\mu(U) & := &\deg U/\mbox{rk}\,U\nonumber\eeqn\noin
is called the \emph{slope} of $U$.\end{definition}

Clearly, if $E$ is stable as a vector bundle  ---  meaning that all of its subbundles satisfy \eqref{Stable}  ---  then for any Higgs field $\phi\in H^0(X;\End\,E\tensor K^*)$, the pair $(E,\phi)$ is also stable.\\

\begin{remark}  An important property of stable co-Higgs bundles is that they are \emph{simple}: if $(E,\phi)$ is stable, then every endomorphism of $E$ that commutes with $\phi$ is a multiple of the identity.  A proof can be quickly adapted from the analogous result for stable vector bundles; see for instance \cite{NS:65}.
\end{remark}

If $(E,\phi)$ is semistable but not stable, $E$ has a proper subbundle $U$ for which $(U,\phi)$ is stable.  It follows that $(E/U,\phi)$ is semistable.  This process, which terminates eventually, gives us a \emph{Jordan--H\"older filtration} of $E$:
\begin{equation*} 0=E_0\subset\cdots\subset E_m=E\end{equation*}
for some $m$, where $(E_j,\phi)$ is semistable for $1\leq i\leq m-1$, and where $(E_j/E_{j-1},\phi)$ is stable and $\mu(E_j/E_{j-1})=\mu(E)$ for $1\leq j\leq m$.  While this filtration is not unique, the isomorphism class of the following object is:\beqn\mbox{gr}(E,\phi) & := & \Plus_{j=1}^m(E_j/E_{j-1},\phi).\nonumber\eeqn\noin  This object is called the \emph{associated graded object} of $(E,\phi)$.  Then, two semistable pairs $(E,\phi)$ and $(E',\phi')$ are said to be \emph{$S$-equivalent} whenever $\mbox{gr}(E,\phi)\cong\mbox{gr}(E',\phi')$.  If a pair is strictly stable, then the underlying bundle has the trivial Jordan--H\"older filtration consisting of itself and the zero bundle, and so the isomorphism class of the graded object is nothing more than the isomorphism class of the original pair.  

For an arbitrary line bundle $L$ in place of $K^*$, the above notions of isomorphism, semistability, and $S$-equivalence are defined identically.

\section{Higher genus}

Stable co-Higgs bundles with sufficiently interesting Higgs fields occur only on the projective line.  To see this, suppose that $X$ has genus $g>1$ and that $(E,\phi)$ is a stable co-Higgs bundle on $X$.  The canonical line bundle $K$ has $g$ 
global sections: choose one, say, $s$.  Taking the product $s\phi$ contracts $K$ with $K^*$; that is, $s\phi$ is an endomorphism of $E$.  But $s\phi$ and $\phi$ commute, and so $s\phi$ must be a multiple of the identity, by the ``simple'' property of stability.  Because $\deg K=2g-2>1$, $s$ vanishes somewhere, and so $\phi$ must vanish everywhere.  In other words, a stable co-Higgs bundle on $X$ with $g>1$ is nothing more than a stable vector bundle.

When $g=1$, co-Higgs bundles are Higgs bundles.

This leaves only the projective line.  We will see that stable co-Higgs bundles with nonzero Higgs fields are plentiful here.  This is in contrast to Higgs bundles, which are never stable on $\PP^1$.  Co-Higgs bundles, therefore, are an extension of the theory of Higgs bundles to genus 0.

\section{Nitsure's moduli space}

For the existence and features of the moduli space we rely on \cite{NN:91}, in which Nitsure constructs a quasiprojective variety that is a coarse moduli space for $S$-equivalence classes of semistable $L$\emph{-pairs} of rank $r$ on a curve $X$.  Here, $L$ is a sufficiently-ample line bundle and ``$L$-pair'' means a pair $(E,\phi)$ in which $E$ is a rank-$r$ vector bundle and $\phi\in H^0(X;\End E\tensor L)$.  The construction uses geometric invariant theory, and the stability condition is the one defined previously.  For $X=\PP^1$ and $L=\CO(2)$, we have the moduli space of semistable co-Higgs bundles on the projective line.  We use $\CM(r)$ to signify this space; $\CM(r,d)$, the locus in $\CM(r)$ consisting of degree-$d$ co-Higgs bundles.  When $r$ and $d$ are coprime, $\CM(r,d)$ is smooth
 and every point is strictly stable.

For $r=2$, we need only describe the loci $\CM(2,-1)$ and $\CM(2,0)$, as we can recover co-Higgs bundles of other degrees by tensoring the elements of these two spaces by $\CO(\pm1)^{\tensor n}$ for an appropriate $n$.  In \cite{NN:91} Nitsure calculates the dimension of $\CM(r)$ to be $2r^2+1$, and so $\CM(2)$ is 9-dimensional.  (He proves that the dimension is independent of $d$.)  For a simplification, we consider only trace-free Higgs fields.  The map
\beqn\CM(2) & \rightarrow & H^0(\PP^1;\CO(2))\times\CM_0(2)\nonumber\eeqn\noin 
defined by\beqn(E,\phi) & \mapsto & \left(\mbox{Tr}\,\phi,\left(E,\phi-\frac{1}{2}\mbox{Tr}\,\phi\right)\right),\nonumber\eeqn\noin where $\CM_0(2)$ denotes the 6-dimensional trace-free part of the moduli space, is an isomorphism.  As $\mbox{Tr}\,\phi$ is a Higgs field for a line bundle, the factorization can be thought of as $\CM(2)\cong\CM(1)\times\CM_0(2)$, where the first factor is the space of co-Higgs line bundles of some fixed degree.  The piece of the moduli space that we do not already understand is $\CM_0(2)$, and so there is no generality lost in restricting attention to it.

\section{Hitchin morphism and spectral curves}

Consider the \emph{Hitchin map} $h:\CM(r)\rightarrow\Plus_{k=1}^rH^0(\PP^1;\CO(2k))$ given by $(E,\phi)\mapsto\mbox{char}\,\phi$, where $\mbox{char}\,\phi$ is the characteristic polynomial of $\phi$.  Since $\mbox{char}\,\phi$ is invariant under conjugation, this map is well-defined on equivalence classes.  Nitsure proves in \cite{NN:91} that $h$ is proper.  In particular, pre-images of points are compact.  Therefore, the fibres of $h$ are compact.

Let $\p=(\p_1,\dots,\p_r)\in\Plus_{k=1}^rH^0(\PP^1;\CO(2k))$ be a generic section.  It follows from more general arguments in \cite{BNR:89USCAN} and \cite{DM:96} that the fibre $h^{-1}(\p)$ is isomorphic to the Jacobian of a \emph{spectral curve} embedded as a smooth subvariety $X_\p$ of the total space of $\CO(2)$.  The correspondence works like this:
{\renewcommand{\theenumi}{\alph{enumi}}
\bitem\item\loclab If $\pi$ is the projection to $\PP^1$ of the total space of $\CO(2)$, then the restriction $\pi_\p:X_\p\rightarrow\PP^1$ is an $r:1$ covering map. 
\item\loclab If $y$ is the coordinate on the total space of $\CO(2)$ and $\eta$ is the tautological section of the pullback of $\CO(2)$ to its own total space, then the equation of $X_\p$ is $\eta^r(y)=\p_1(\pi(y))\eta^{r-1}(y)+\cdots+\p_r(\pi(y))$.
\item\loclab The direct image of a line bundle $L$ on a generic $X_\p$ is a rank-$r$ vector bundle $(\pi_\p)_*L=E$ on $\PP^1$.
\item\loclab The pushforward of the multiplication map $L\rightarrow\eta L$ is a Higgs field $\phi$ for $E$, with characteristic polynomial $\p$.\eitem
We admit that we are abusing language, by referring to $\p$ as the characteristic polynomial when it is the tuple of characteristic coefficients. }

The spectral curve ramifies at finitely-many points, which are the $z\in\PP^1$ for which $\phi_z$ has repeated eigenvalues.  The generic characteristic polynomial $\p$ is irreducible, and so its $X_\p$ is an irreducible curve.

In the case of rank $r=2$ and $\phi$ trace-free, the characteristic polynomial is a monic polynomial of degree 2 in $\eta$ with no linear term, and with a section of $\CO(4)$ for the coefficient of $\eta^0$.  This section vanishes at 4 generically distinct points in $\PP^1$, which are the ramification points of the double cover $X_\p\rightarrow\PP^1$.  By the Riemann--Hurwitz formula, $X_\p$ is an elliptic curve, whose Jacobian is another elliptic curve.  Therefore, the map $h$ on $\CM_0(2)$ is a fibration of generically nonsingular elliptic curves over the 5-dimensional affine space of determinants. 

Because the generic $X_\p$ is irreducible, a co-Higgs bundle $(E,\phi)$ coming from a line bundle on $X_\p$ has no $\phi$-invariant subbundles whatsoever, and therefore is stable.  Stability limits the underlying vector bundles that can be obtained from spectral line bundles.  In the next section, we address this.

\section{Stable Grothendieck numbers}

According to the classical Birkhoff--Grothendieck theorem, if $E$ is a rank-$r$ holomorphic vector bundle on $\PP^1$, then 
\[E\cong\CO(m_1)\plus\CO(m_2)\plus\cdots\CO(m_r)\] 
for integers $m_1,m_2,\dots,m_r$ that are unique up to permutation. We find necessary and sufficient conditions on the Grothendieck numbers $m_i$ for the existence of semistable Higgs fields.

\begin{theorem}\label{MainResult}  Let $E=\CO(m_1)\plus\CO(m_2)\plus\cdots\plus\CO(m_r)$ be a holomorphic vector bundle of rank $r>1$ on $\PP^1$.  If the line bundles are ordered so that $m_1\geq m_2\geq\cdots\geq m_r$,  then $E$ admits a semistable $\phi\in H^0(\PP^1;\mbox{\emph{End}}\,E\tensor\CO(2))$ if and only if $m_i\leq m_{i+1}+2$ for all $1\leq i\leq r-1$.  The generic $\phi$ leaves invariant no subbundle of $E$ whatsoever; therefore, the generic $\phi$ is stable trivially.\end{theorem}

\begin{proof}  We begin with the \emph{only if} direction, for which we proceed by induction on successive extensions of balanced bundles by each other.  (A rank-$r$ \emph{balanced vector bundle} over $\PP^1$ splits into $r$ copies of a single line bundle.)  To arrive at these bundles, we filter the decomposition of $E$ by its repeated Grothendieck numbers.  That is, if the first $d_1$ ordered Grothendieck numbers are $m_1=\cdots=m_{d_1}=a_1$, then we write $E_1$ for the balanced vector bundle $\Plus^{d_1}\CO(a_1)$.  If the next $d_2$ numbers are all equal to the same number, say $a_2$, then we set $E_2:=\Plus^{d_2}\CO(a_2)$; and so on.  Then, $E=\Plus_{i=1}^kE_i=\Plus_{i=1}^k\left(\Plus^{d_i}\CO(a_i)\right)$, where $d_1+\cdots+d_k=r$ and $a_1>\cdots>a_k$.

Begin with the sequence
\beqn E_1\stackrel{\phi}{\rightarrow}E\tensor\CO(2)\stackrel{p}{\rightarrow}(E_2\plus\cdots\plus E_k)\tensor\CO(2).\nonumber\eeqn\noin  The composition of $\phi$ with the quotient map $p$ is a  section of 
\[E_1^*\tensor(E/E_1)\tensor\CO(2),\] 
and so has components in $\CO(-a_1+a_j+2)$, for each of $j=2,3,\dots,k$.  If $a_1>a_2+2$, then $a_1>a_j+2$ for $j=2,3,\dots,k$ and
\beqn H^0(\PP^1;\CO(-a_1+a_2+2))=\cdots=H^0(\PP^1;\CO(-a_1+a_k+2))=0.\nonumber\eeqn\noin  Therefore, $p\circ\phi$ is the zero map.  It follows that $E_1$ is $\phi$-invariant, and since $d_1+\cdots+d_k=r$ and $a_1>a_2>\cdots>a_k$, we have
\beqn\frac{\deg E_1}{\mbox{rk}\,E_1}=\frac{d_1a_1}{d_1}=a_1=\frac{a_1(d_1+\cdots+d_k)}{r}\!>\!\frac{d_1a_1+d_2a_2+\cdots d_ka_k}{r}=\frac{\deg E}{\mbox{rk}\,E}.\nonumber\eeqn

Because $(E,\phi)$ is semistable, such a subbundle of $E$ cannot exist.  In light of the contradiction, we must have $a_1\leq a_2+2$, and so
\beqn m_1=\cdots=m_{d_1} & \leq  & m_{d_1+1}+2=\cdots=m_{d_1+d_2}+2.\nonumber\eeqn 
Assume now that
\beqn a_2 & \leq & a_3+2\nonumber\\ &\;\; &\vdots  \nonumber\\a_{j-1} & \leq & a_j+2,\nonumber\eeqn\noin 
and examine the sequence
\beqn E_1\plus E_2\plus\cdots\plus E_j\stackrel{\phi}{\rightarrow}E\tensor\CO(2)\stackrel{p}{\rightarrow}(E_{j+1}\plus\cdots\plus E_k)\tensor\CO(2)\nonumber\eeqn\noin 
in which we abuse notation and re-use $p$ for the quotient of $E$ by $E_1\plus\dots\plus E_j$.  We assume that $a_j>a_{j+1}+2$.  Because of the induction hypothesis, we have that $a_i\geq a_j>a_u+2$ for each $i\leq j$ and each $u>j$.  Therefore, $-a_i+a_u+2<0$, and the images of the balanced bundles $E_i$, $i\leq j$, are zero under the composition of $\phi$ and $p$.  Hence, $E_1\plus\cdots\plus E_j$ is $\phi$-invariant and its slope exceeds that of $E$.  The induction is complete.

\begin{remark} The validity of the argument above is not exclusive to $X=\PP^1$: $X$ could be projective space $\PP^n$ of any dimension, so long as we are considering fully decomposable bundles.  In that case, the result would say that semistable Higgs fields exist only if $m_i\leq m_{i+1}+s$, where $s$ is the largest integer such that $T(-s)$ has sections.
\end{remark}

Conversely, suppose that $m_i\leq m_{i+1}+2$ for each $i=1,\dots,r-1$.  Our strategy is to find a particular Higgs field $\phi$ under which no subbundle of $E$ is invariant, meaning that $(E,\phi)$ is trivially stable.  Because of the decomposition of $E$ into a sum of line bundles $\CO(m_i)$, the Higgs field can be realized as an $r\times r$ matrix whose $(i,j)$-th entry takes values in the line bundle $\CO(-m_j+m_i+2)$.  The subdiagonal elements are sections of $\CO(-m_{i-1}+m_{i}+2)\cong\CO(p_i)$ for $i=2,\dots,r$, where each $p_i$ is one of 0, 1, or 2.  Into each of these positions, we enter a `$1$', which represents the section of $\CO(p_i)$ that is $1$ on $\PP^1-\set{\infty}$ and is $1/z^{p_i}$ on $\PP^1-\set{0}$.  The $(1,r)$-th entry is a section of $\CO(-m_r+m_1+2)$, which is of degree 2 or more.  There, we insert $z$.  For all other entries, we insert the zero section of the corresponding line bundle:\bcent\begin{math}\phi(z)=\left(\begin{array}{cccccc}0 & 0 & \cdots & 0 & 0 & z\\1 & 0 & \cdots & 0 & 0 & 0\\ \vdots & \vdots & \ddots & \vdots & \vdots & \vdots \\ 0 & 0 & \cdots & 1 & 0 & 0\\ 0 & 0 & 0 & \cdots & 1 & 0\end{array}\right).\end{math}\ecent
   \noin Over $\PP^1-\set{\infty}$, the characteristic polynomial of $\phi$ is $(-1)^{r-1}z+y^r$, which is irreducible in $\CC[y][z]$.

Because the characteristic polynomial does not split, $\phi$ has no proper eigen-subbundles in $E$; that is, $E$ has no $\phi$-invariant subbundles.  As irreducibility is an open condition, the genericity follows immediately: there is a Zariski open subset of $H^0(\PP^1;\End E\tensor\CO(2))$ whose elements leave invariant no subbundles of $E$ whatsoever.\end{proof}

For the case of rank $r=2$, Theorem \ref{MainResult} tells us that if $E$ has degree 0, then $E$ admits semistable Higgs fields if and only if $E\cong\CO\plus\CO$ or $E\cong\CO(1)\plus\CO(-1)$.  On the other hand, if $E$ has degree $-1$, there is only one choice: $E\cong\CO\plus\CO(-1)$.

\section{Odd degree}

We examine $\CM_0(2,-1)$, where the underlying bundle of every co-Higgs bundle is isomorphic to $E=\CO\plus\CO(-1)$.  Since $E$ has non-integer slope, every semistable Higgs field for $E$ is stable.  Every Higgs field for $E$ is of the form\bcent\begin{math}\phi=\left(\begin{array}{cr}a & b\\c & -a\end{array}\right),\end{math}\ecent\noin where $a$, $b$, and $c$ are sections of $\CO(2)$, $\CO(3)$, and $\CO(1)$, respectively.  The stability of $\phi$ means that $c$ is not identically zero: because $\mu(E)=-1/2$, $\phi$ cannot leave the trivial sub-line bundle $\CO$ invariant.  Accordingly, $c$ has a unique zero $z_0\in\PP^1$.

It is possible to provide a global description of the odd-degree moduli space as a universal elliptic curve.  Let $\pi:M\rightarrow\PP^1$ stand for the two-dimensional total space of $\CO(2)$.  We claim that we can assign uniquely to each stable $\phi$ a point in the 6-dimensional space $\CS$ defined by
\beqn\set{(y,\p)\in M\times H^0(\PP^1;\CO(4)) \; : \; \eta^2(y)=\p(\pi(y))}.\nonumber\eeqn\noin  
That $\CS$ is a smooth subvariety of the 7-dimensional space $M\times H^0(\PP^1;\CO(4))$ can be seen as follows.  Over the subset $U_0$ of $\PP^1$ where the coordinate $z$ is not $\infty$, we have
\beqn\CS & = & \set{(z,y,a_0,a_1,a_2,a_3,a_4) : y^2=a_0+a_1z+a_2z^2+a_3z^3+a_4z^4},\label{Model}\eeqn
 with $(z,y)$ as coordinates on $M$.  If $\tilde z=1/z$ and $\tilde y=y/z^2$, then $(\tilde z,\tilde y)$ give coordinates on $M$ over $U_1=\PP^1-\set{0}$.  There, $\CS$ is given by 
\beqn\tilde y^2 & = & a_4+a_3\tilde z+\cdots+a_0\tilde z^4.\nonumber\eeqn
Since $\di f/\di a_0\neq 0$ on $M|_{U_0}\times\CC^5$ and $\di\tilde f/\di a_4\neq0$ on $M|_{U_1}\times\CC^5$, where 
\beqn
f(z,y,a_0,\dots,a_4)&=&y^2-a_0-a_1z-\cdots-a_4z^4,\nonumber\\ 
\tilde f(\tilde z,\tilde y,a_0,\dots,a_4)&=&\tilde y^2-a_4-a_3\tilde z-\cdots-a_0\tilde z^4,\nonumber
\eeqn\noin
the variety $\CS$ is in fact smooth as a subvariety.

We will define an isomorphism from $\CM_0(2,-1)$ onto $\CS$ by sending $\phi$ to $(z_0,a(z_0),-\det\phi)$, with $z_0$ and $a$ as above.  Since $a$ is a section of $\CO(2)$, $(z_0,a(z_0))$ is a point in $M$.  The point is determined uniquely by the conjugacy class of $\phi$, for if 
\[
\psi=\left(\begin{array}{cc}d & e\\0 & f\end{array}\right)
\]
is an automorphism of $E=\CO\plus\CO(-1)$, in which case $e$ is a section of $\CO(1)$ and $d,f\in\CC^*$, then the Higgs field transforms as
\[
\phi'=\psi\phi\psi^{-1}=\left(\begin{array}{cc}a+d^{-1}ec & -f^{-1}(2ea-bd+e^2fc)\\d^{-1}fc & -a-d^{-1}ec\end{array}\right).
\]
Because $(a+d^{-1}ec)(z_0)=a(z_0)$, the image of $\phi$ in the variety $\CS$ remains unchanged by $\phi\rightarrow\phi'$.  Furthermore, we have $(a(z_0))^2=-\left.\det\phi\right|_{z=z_0}$, and therefore $(z_0,a(z_0),-\det\phi)$ is a point in $\CS$.

Now we start with a point $(z_0,y_0,a_0,a_1,a_2,a_3,a_4)\in M\times\CC^5$.  To be in $\CS$, the point must have $y_0^2=a_0+a_1z_0+\cdots+a_4z_0^4$.  There are two choices of $y_0$, corresponding to the two square roots of $a_0+a_1z_0+\cdots+a_4z_0^4$, unless  $a_0+a_1z_0+\cdots+a_4z_0^4=0$, in which case the point in $\CS$ is $(z_0,0,0,0,0,0,0)$.  Let us assume for the moment that $z_0$ is such that $a_0+a_1z_0+\cdots+a_4z_0^4\neq0$.  The two corresponding points in $\CS$ are
\beqn\left(z_0,\sqrt{a_0+a_1z_0+a_2z_0^2+a_3z_0^3+a_4z_0^4\;},a_0,a_1,a_2,a_3,a_4\right)\nonumber\eeqn\noin 
and
\beqn \left(z_0,-\sqrt{a_0+a_1z_0+a_2z_0^2+a_3z_0^3+a_4z_0^4\;},a_0,a_1,a_2,a_3,a_4\right).\nonumber\eeqn\noin  
Consider the first of the two points.  Its pre-image in $\CM_0(2,-1)$ is a stable Higgs field
\bcent\begin{math}\phi=\left(\begin{array}{cr}a & b\\c & -a\end{array}\right)\end{math}\ecent\noin 
for which $z_0$ is the unique point in $\PP^1$ at which $c$ vanishes, 
\beqn\det\phi=-a_0-a_1z-\cdots-a_4z^4,\eeqn
and $a(z_0)=y_0$.  A representative Higgs field has
\beqn a & = & \sqrt{a_0+a_1z_0+a_2z_0^2+z_3z_0^3+a_4z_0^4\;},\nonumber\\b(z) & = & a_1+a_2z_0+a_3z_0^2+a_4z_0^3+(a_2+a_3z_0+a_4z_0^2)z\notag
\\&  &\quad\;+(a_3+a_4z_0)z^2+a_4z^3,\nonumber\\c(z) & = & z-z_0.\nonumber\eeqn\noin  If we use $a=-\sqrt{a_0+a_1z_0+a_2z_0^2+z_3z_0^3+a_4z_0^4\;}$ instead, then we get a Higgs field for the other point in $\CS$. 

For convenience, choose a coordinate $z$ that vanishes at $z_0$.  Then, the two points in $\CS$ are $\left(z_0,\sqrt{a_0\,},a_0,a_1,a_2,a_3,a_4\right)$ and $\left(z_0,-\sqrt{a_0\,},a_0,a_1,a_2,a_3,a_4\right)$, and their respective Higgs fields become\bcent\begin{math}\phi_+(z)=\left(\begin{array}{cc}\sqrt{a_0} & a_1+a_2z+a_3z^2+a_4z^3\\ z & -\sqrt{a_0}\end{array}\right)\end{math}\ecent\noin and\bcent\begin{math}\phi_-(z)=\left(\begin{array}{cc}-\sqrt{a_0} & a_1+a_2z+a_3z^2+a_4z^3\\ z & \sqrt{a_0}\end{array}\right)\end{math}.\ecent  The two points coincide with each other, and $\phi_+=\phi_-$, when $a_0=0$.  Having $a_0=0$ is equivalent to the spectral curve ramifying above $z_0$, because $a_0=0$ means that the characteristic equation of $\phi_{\pm}$ is 
\[y^2=-z(a_1+a_2z+a_3z^2+a_4z^3),\] 
and so $y^2=0$ at $z_0$.

Since $\phi_+$ and $\phi_-$ correspond to distinct points in $\CS$ whenever $z_0$ is not a ramification point of their corresponding spectral curve, there can be no automorphism of $E=\CO\plus\CO(-1)$ that takes $\phi_+$ to $\phi_-$, unless $a_0=0$.  This is easy to verify. Suppose that there exists a $\psi\in H^0(\mbox{Aut}\,E)$, say
\bcent\begin{math}\psi=\left(\begin{array}{cc}d & e\\0 & f\end{array}\right)\end{math}\ecent\noin 
with $d,f\in\CC^*$ and $e$ a section of $\CO(1)$, such that $\psi\phi_+\psi^{-1}=\phi_-$.  The matrix $\psi\phi_+\psi^{-1}$ is 
\bcent\begin{math}\displaystyle\frac{1}{df}\left(\begin{array}{cc}df\sqrt{a_0}+efz & -2de\sqrt{a_0}+d^2\tilde b(z)-e^2z\\f^2z & -df\sqrt{a_0}-efz\end{array}\right),\end{math}\ecent\noin in which $\tilde b(z)=a_1+a_2z+a_3z^2+a_4z^3$. Equality with $\phi_-$ requires $f=d$ and $2\sqrt{a_0}=-\displaystyle\frac{e}{d}z$.  Since $-e/d$ is a section of $\CO(1)$ we can write it as $lz+m$ for some $l,m\in\CC$, and so the condition becomes $2\sqrt{a_0}=lz^2+mz$.  This can only be satisfied when $a_0=0$ (and $l=m=0$).  

We can frame this discussion by appealing to the spectral viewpoint.  Consider a generic spectral curve, which is a smooth curve of genus 1.  According to Grothendieck--Riemann--Roch, to get $E=\CO\plus\CO(-1)$ on $\PP^1$, we need a degree-1 line bundle $L$ on the spectral curve.   The ordinary Riemann--Roch theorem tells us that $L$ has a one-dimensional space of global holomorphic sections, and so all of these sections must vanish at a single point.  Using the coordinates on $M$, this point is either $(z_0,\sqrt{a(z_0)}\,)$ or $(z_0,-\sqrt{a(z_0)}\,)$.  Whether we have $\phi_+$ or $\phi_-$ depends on which sheet of the double cover contains the point at which the sections of $L$ vanish.  The covering map for the spectral curve projects $(z_0,\sqrt{a(z_0)}\,)$ and $(z_0,-\sqrt{a(z_0)}\,)$ onto $z_0$, the point in $\PP^1$ at which the $\CO(1)$-components of $\phi_+$ and $\phi_-$ vanish.  If the vanishing point of the global sections of $L$ is a point where the two sheets coincide, then we get a single stable Higgs field $\phi_+=\phi_-$.

Our construction of $\phi_\pm$ and our argument regarding automorphisms of $E$ are independent of whether the spectral curve is singular or nonsingular, and so our isomorphism $\CM_0(2,-1)\cong\CS$ holds globally.

\section{Even degree}

The moduli space $\CM_0(2,0)$ does not yield
such an explicit description; however, we can still say something about the fibres of the Hitchin map.

Recall that Theorem \ref{MainResult} allows for two choices of underlying bundle: $E^{\;\;\,1}_{-1}:=\CO(1)\plus\CO(-1)$ or the trivial rank-2 bundle $E_0:=\CO\plus\CO$, the latter of which is the generic splitting type.  
If a pair $(E^{\;\;\,1}_{-1},\phi)$ is \emph{not} unstable, then it is strictly stable: every sub-line bundle of degree 0 is contained in $\CO(1)$, and is therefore $\phi$-invariant if and only if $\CO(1)$ is $\phi$-invariant. 
On the other hand, $E_0$ admits semistable but not stable Higgs fields $\phi$: these are the upper-triangular Higgs fields, in which the three matrix coefficients in the polynomial $\phi(z)=A_0+A_1z+A_2z^2$ admit a common eigenvector.  
 The $S$-equivalence class of such a $\phi$ is represented by the graded object\bcent\begin{math}\mbox{gr}(\phi)=\left(\begin{array}{cc}a & 0\\0 & -a\end{array}\right),\end{math}\ecent\noin for some $a\in H^0(\PP^1;\CO(2))$.  Consequently, every point in a generic fibre of the Hitchin map is strictly stable, because $\p=-a^2$ is a reducible spectral curve, whereas the generic spectral curve is irreducible.
One example of a non-generic fibre is the nilpotent cone over $\p=0$: in addition to stable Higgs fields it also contains the zero Higgs field for $E_0$, which is semistable but not stable.

To study Higgs fields for $E^{\;\;\,1}_{-1}$, we define a section of the Hitchin map $h:\CM_0(2,0)\rightarrow H^0(\PP^1;\CO(4))$ in the following way: to each $\p\in H^0(\PP^1;\CO(4))$, we assign the Higgs field\bcent\begin{math}Q(\p)=\left(\begin{array}{cc}0 & -\p\\1 & 0\end{array}\right)\end{math}\ecent\noin for $E^{\;\;\,1}_{-1}$, with the symbol 0 denoting the zero section of $\CO(2)$, and where 1 is unity.  This section is the genus-0 analogue of Hitchin's model of Teichm\"uller space \cite{NJH:86}, but with our $\p$ replacing the quadratic differential in his model.

\begin{proposition}  The section $Q$ is the locus in $\CM_0(2,0)$ of stable co-Higgs bundles with \hbox{underlying} bundle isomorphic to $E^{\;\;\,1}_{-1}=\CO(1)\plus\CO(-1)$.\end{proposition}

\begin{proof}If \bcent\begin{math}\phi=\left(\begin{array}{cr}a & b\\c & -a\end{array}\right)\end{math}\ecent\noin is a stable Higgs field for $E^{\;\;\,1}_{-1}$, then $a$ is a section of $\CO(2)$, $b$ is a section of $\CO(4)$, and $c$ is a constant.  Stability implies that $c\neq0$.  To study the orbit of $\phi$ under automorphisms of $E^{\;\;\,1}_{-1}$, we take\bcent\begin{math}\psi=\left(\begin{array}{cc}1 & d\\0 & e\end{array}\right)\rlap,\end{math}\ecent\noin in which $d$ is a section of $\CO(2)$ and $e\in\CC^*$.  The transformed Higgs field is\bcent\begin{math}\phi'=\psi\phi\psi^{-1}=\left(\begin{array}{cc}a+dc & -2de^{-1}a+e^{-1}b-d^2e^{-1}c\\ec & -a-dc\end{array}\right).\end{math}\ecent\noin  Taking the automorphism $\psi$ with $e=c^{-1}$, $d=-ac^{-1}$, we get\bcent\begin{math}\phi'=\psi\phi\psi^{-1}=\left(\begin{array}{cc}0 & a^2+bc\\1 & 0\end{array}\right).\end{math}\ecent\noin  In other words, the conjugacy class of a trace-free Higgs field acting on $E^{\;\;\,1}_{-1}$ is determined by a unique $\p=a^2+bc=-\det\phi\in H^0(\PP^1;\CO(4))$.\end{proof}
~~

Consider a generic spectral curve $X_\p$, which again is a smooth curve of genus 1.  Grothendieck--Riemann--Roch tells us the following: for the direct image of a line bundle $L$ on $X_\p$ to be a rank-2 vector bundle of degree 0 on $\PP^1$, then we must have $\deg L=2$.  On $\PP^1$, twisting $E_0$ by $\CO(-1)$ gives $\CO(-1)\plus\CO(-1)$, which has no global sections.  On the other hand, twisting $E^{\;\;\,1}_{-1}$ by $\CO(-1)$ gives $\CO\plus\CO(-2)$, which still has a global section.  Because the direct image functor preserves the number of global sections, this is the same as asking whether or not $L\tensor\pi_\p^*\CO_{\PP^1}(-1)$ has global sections.  The twisted line bundle $L\tensor\pi_\p^*\CO_{\PP^1}(-1)$ has degree $\deg L+(-1)\deg\pi_\p=2-2=0$.  The only line bundle of degree 0 on $X_\p$ with a global section is the trivial line bundle $\CO_{X_\p}$.  Therefore, pushing down $\CO_{X_\p}\tensor\pi_\p^*\CO_{\PP^1}(1)$ produces the co-Higgs bundle $(E^{\;\;\,1}_{-1},Q(\p))$, while pushing down any other line bundle of degree 2 gives a Higgs field for $E_0$.

\section{Betti numbers and holomorphic chains}\label{SecBetti}

In this section, we reincorporate the trace of $\phi$; that is, we consider the full moduli space $\CM(r,d)$ of stable rank-$r$ and degree-$d$ co-Higgs bundles on $\PP^1$.

As with the conventional Higgs bundle moduli space, $\CM(r,d)$ enjoys a circle action, $(E,\phi)\mapsto(E,e^{i\theta}\phi)$, which induces a localization of the Poincar\'e series of $\CM(r,d)$ whenever $\gcd(r,d)=1$.  This localization originates in Morse--Bott theory and is developed in \cite{NJH:86,PBG:94,PBG:95,TH:98,HT:03} for the case of the Hitchin system.  All of the arguments carry over to co-Higgs bundles without modification.  (The Morse--Bott function, defined to be a scalar multiple of the norm squared of the Higgs field using the natural K\"ahler metric, is a proper moment map for the action and is perfect and nondegenerate, as discussed in Proposition 7.1 and Theorem 7.6 of \cite{NJH:86}.  
We will not need to interact with the function directly.)  Before we state the main features of the theory, we need the following notion: 
if $k$ is a nonnegative integer and $(U_1,\dots,U_n)$ is an ordered $n$-tuple of vector bundles such that $E=\Plus U_i$, then an element
$\psi\in H^0(\End E)$ is said to act with \emph{weight} $k$ on $(U_1,\dots,U_n)$ if $\psi(U_i)\subseteq U_{i+k}$.  (If $i+k>n$, then $\psi(U_i)=0$.)  This notion extends to twisted morphisms as well, that is, when $\psi\in H^0(\End E\tensor L)$ for some line bundle $L$.

Now, the main features of Morse theory for co-Higgs bundles, as adapted from Higgs bundles, are:\bitem

\item The downward gradient flow of the Morse--Bott function is coincident with the nilpotent cone, and the moduli space deformation retracts onto the cone (\S4.4 of \cite{TH:98}).

\item A fixed point of the circle action is a co-Higgs bundle with a special form: a \emph{holomorphic chain}.  Such an object is a $(2n-1)$-tuple $(U_1,\dots,U_n;\phi_1,\dots,\phi_{n-1})$ for some $n\leq r$, in which each $U_i$ is a holomorphic vector bundle on $\PP^1$ and each $\phi_i$ is a holomorphic map $U_i\rightarrow U_{i+1}\tensor\CO(2)$. (We refer to \S7 of \cite{NJH:86} for $r=2$, and to Lemma 2 of \cite{TSCH:91I} and p.18 of \cite{PBG:95} for higher rank.)  For the case of ordinary Higgs bundles, where the twist is by the canonical line bundle, these objects are complex variations of Hodge structure \cite{TSCH:91I}.  The term ``holomorphic chain'', which accommodates more general twisting, originates in \cite{ACGP:01}.

The \emph{total rank} of a chain is $\sum\rk U_i$. Its \emph{total degree} is $\sum\deg U_i$.  The \emph{type} of a chain is the vector $(\rk U_1,\dots,\rk U_n)$ and the \emph{degree vector} is $(\deg U_1,\dots,\deg U_n)$.  By taking $E=\Plus U_i$ and writing down a block matrix $\phi$ with sub-diagonal blocks $[\phi]_{i+1,i}=\phi_i$ and zero blocks elsewhere, we get a nilpotent co-Higgs bundle.  The Higgs field is an element of $H^0(\End E\tensor\CO(2))$ acting with weight 1 on $(U_1,\dots,U_n)$.  We define a chain to be \emph{(semi)stable} when its associated Higgs bundle is (semi)stable.  Since $\gcd(r,d)=1$, all of the chains we shall consider are strictly stable.\enlargethispage{\baselineskip}

\item The \emph{Morse index} at a fixed point is the number of negative eigenvalues of the Hessian of the Morse--Bott function at the fixed point. Let $(E,\phi)$ be any fixed point, with decomposition $E=\Plus U_i$.  We denote its Morse index by $\beta((E,\phi))$.  After a calculation involving the Hessian, Gothen shows (p.19 of \cite{PBG:95}) that $\beta((E,\phi))$ is a sum of two integers $\beta^{1,0}$ and $\beta^{0,1}$, where $\beta^{1,0}$ is the real dimension of the subspace in
\beqn\quad\quad\quad\quad\;\frac{H^0(\End E\tensor\CO(2))}{\img H^0(\End E)\stackrel{[-,\phi]}{\longrightarrow}H^0(\End E\tensor\CO(2))}\label{Top}\eeqn\noin 
consisting of the elements acting with weight $\geq2$ on $(U_1,\dots,U_n)$, and $\beta^{0,1}$ is the real dimension of the subspace of
\beqn\quad\quad\quad\quad\;\ker H^1(\End E)\stackrel{[-,\phi]}{\longrightarrow}H^1(\End E\tensor\CO(2))\label{Bottom}\eeqn\noin
consisting of elements acting with weight $\geq 1$ on $(U_1,\dots,U_n)$.

It is an immediate consequence of the stable $\Rightarrow$ simple property and Serre duality that the map in the denominator of \eqref{Top} is injective, and the map in \eqref{Bottom} is surjective.  This means that the Morse index is given by the more compact formula 
\beqn\quad\quad\;\beta((E,\phi))  &  = & {\dim_\RR H^0_{\geq2}(\End E\tensor\CO(2))-\dim_\RR H^0_{\geq1}(\End E))\nonumber}\\ &  &\quad\; +\dim_\RR H^1_{\geq1}(\End E)-\dim_\RR H^1_{\geq2}(\End E\tensor\CO(2)),\nonumber\eeqn\noin 
where the subscripts $\geq k$ refer to which weight subspace is being isolated.  Notice that there are two Riemann--Roch identities interlaced by this formula.  Taking advantage of this allows us to further reduce the formula to
\beqn & & \beta((E,\phi)) \;\;=\;\;  4\delta^n_2\sum_{i=1}^{n-2}\sum_{j=i+2}^nr_ir_j-2\delta^n_1\sum_{i=1}^{n-1}(\deg(U_i^*U_{i+1})+r_ir_{i+1}),\label{TheFormula}\eeqn\noin 
where $r_i=\rk U_i$ and $\delta^n_j=1$ if $n>j$ and $0$ otherwise.  (Note that $\deg(U_i^*U_{i+1})=-r_{i+1}\deg U_i+r_i\deg U_{i+1}$, and so the formula for the Morse index at a fixed point depends only on the ranks and degrees of the bundles in the chain.)\\\eitem

For calculation purposes, it is useful to know that the chain type and degree vector are constant on connected components of the fixed point set. (This is Lemma 9.2 in \cite{HT:03}, which is attributed by the authors to Carlos Simpson.)  Therefore, to each component of the fixed point set, we may associate a vector $\mb{r}=(r_1,\dots,r_n)\in\ZZ^n_{>0}$ and a vector $\mb{d}=(d_1,\dots,d_n)\in\ZZ^n$, with $\sum r_i=r$ and $\sum d_i=d$.  Since the Morse index depends only on $\mb{r}$ and $\mb{d}$, as in formula \eqref{TheFormula}, we have that the Morse index is constant on connected components of the fixed point set.

The main tool for our calculation of Betti numbers is the Morse-theoretic localization formula (\S7 of \cite{NJH:86}): the Poincar\'e series of $\CM(r,d)$ is
\beqn \CP(r,d;x) & = & \sum_\CN x^{\beta(\CN)}\CP(\CN;x),\nonumber\eeqn\noin 
where $\CN$ stands for a connected component of the fixed point set of the circle action; $\CP(\CN;x)$, for the Poincar\'e polynomial of $\CN$; and $\beta(\CN)$, for the Morse index of  any point in  $\CN$.

Two different connected components can have the same $\mb{r}$ and the same $\mb{d}$.  This will occur when the set of stable chains with type $\mb{r}$ and degree vector $\mb{d}$ is disconnected inside $\CM(r,d)$.  However, two such components will have the same Morse index, determined by $\mb{r}$ and $\mb{d}$.  Therefore, we can rewrite the localization formula as
\beqn\CP(r,d;x) & = & \sum_{\mb r,\mb d}x^{\beta(\mb{r},\mb{d})}\sum_{i\in I(\mb{r},\mb{d})}\CP(\CN_i;x),\nonumber\eeqn\noin 
in which:
\bitem\item The outer sum is taken over all vectors $\mb{r}=(r_1,\dots,r_n)\in\ZZ^n_{>0}$ and $\mb{d}=(d_1,\dots,d_n)\in\ZZ^n$ with $\sum r_i=r$ and $\sum d_i=d$.

\item The exponent $\beta(\mb{r},\mb{d})$ is the right side of formula \eqref{TheFormula} evaluated at $\mb{r}$ and $\mb{d}$.

\item $I(\mb{r},\mb{d})$ indexes the connected components of the set of chains in $\CM(r,d)$ with type $\mb{r}$ and degree vector $\mb{d}$.

\item $\CN_i$ is a connected component of the set of stable chains of type $\mb r$ and degree vector $\mb d$.
\item If the set of stable chains of type $\mb r$ and $\mb d$ is empty, then we declare its Poincar\'e series to be 0.\\
\eitem

\begin{remark} 
Since the nilpotent cone is a deformation retract of $\CM(r,d)$, the Betti numbers of $\CM(r,d)$ and $\CM_0(r,d)$ will be identical.\\
\end{remark}

If we wish to calculate Betti numbers of $\CM(r,d)$, we need to determine all of the stable chains with total rank $r$ and total degree $d$.  Note that for $r>1$, there are no stable chains of type $(r)$, as these are vector bundles on $\PP^1$ with the zero Higgs field.

For $r=2$, there is only one chain type to consider: $(1,1)$.  For $\CM(2,-1)$ in particular, chains of this type have the form $(\CO(a),\CO(-a-1);\phi_1)$ for some integer $a$, where $\phi_1\in H^0(\CO(a)^*\tensor\CO(-a-1)\tensor\CO(2))=H^0(\CO(-2a+1))$.  If $a>0$, then $\phi_1$ must be zero, since $\CO(-2a+1)$ has no global holomorphic sections.  This means that $\CO(a)$ is an invariant sub-line bundle of positive slope in the associated Higgs bundle, which has slope $-1/2$.  If $a<0$, there is an invariant sub-line bundle of degree $-a-1\geq0$, which is also destabilizing.  Therefore, stability necessitates $a=0$, in which case $\phi_1\in H^0(\CO(1))=\CC^2$.  If $\phi_1=0$, then $\CO$ is a destabilizing sub-line bundle.  If $\phi_1\neq0$, then the only invariant subbundles are those contained in $\CO(-1)$, and so their degrees are strictly less than $-1/2$.  Therefore, the stable chains are precisely those of the form $(\CO,\CO(-1);\phi_1)$ with $\phi_1\neq0\in\CC^2$.

Automorphisms of $\CO\plus\CO(-1)$ preserving the chain structure are para\-me\-trized by $\mbox{Aut}(\CO)\plus\mbox{Aut}(\CO(-1))=\CC^*\plus\CC^*$.  The quotient of $\CC^2\bsh\set{0}$ by either right multiplication by the first summand or left multiplication by the second summand of $\CC^*\plus\CC^*$ gives us a connected fixed point set, isomorphic to $\PP^1$.  Because there is only one component of the fixed point set, this $\PP^1$ is the minimal component, whose Morse index is 0.  Putting this together, we have $\CP(2,-1;x)=1+x^2$.

Since the downward Morse flow and the nilpotent cone are coincident, the cone is therefore isomorphic to $\PP^1$.  This is consistent with our concrete model \eqref{Model}. The nilpotent cone in $\CM(2,-1)$ is the subvariety of $\CS$ consisting of points of the form $(z,0,0,0,0,0,0)$.  Since $z$ is just the coordinate on the base $\PP^1$, the nilpotent cone is a copy of $\PP^1$.  Since the moduli space deformation retracts onto the cone, we can read from our model that the Betti numbers of $\CM(2,-1)$ are those of the 2-sphere.

For chains of type $(1,\dots,1)$, we can generalize the discussion from rank 2 to higher rank and arbitrary degree.  If $(L_1,\dots,L_r;\phi_1,\dots,\phi_{r-1})$ is a chain of type $(1,\dots,1)$, then the stability condition is equivalent to $\phi_i\neq0$ for $1\leq i\leq n-1$ and\\ 
\bcent\bmath\begin{array}{ccc}\deg L_r & < & d/r\\ & & \\\displaystyle\frac{\deg L_{r-1}+\deg L_r}{2} & < & d/r\\& \vdots &\\\displaystyle\frac{\deg L_2+\cdots+\deg L_r}{r-1} & < & d/r.\\&&\end{array}\emath\ecent\noin If one of the maps $\phi_i$ were zero, then $L_1\plus\cdots\plus L_i$ and $L_{i+1}\plus\cdots\plus L_r$ would be subbundles of $E$ that are invariant under the associated Higgs field $\phi$.  It is easy to show that they cannot simultaneously have slopes less than $d/r$.  When every $\phi_i$ is nonzero, the slopes of any remaining invariant subbundles are constrained by the inequalities above. Note that the condition $\phi_i\neq0$ requires that $-\deg L_i+\deg L_{i+1}+2\geq0$ for $1\leq i\leq n-1$.  The set of all chains on $(L_1,\dots,L_r)$ is an iterated bundle of projective spaces.  If $L_i=\CO(d_i)$ for each $i$, then the Poincar\'e series of this iterated bundle is equal to the Poincar\'e series of the product $\PP^{-d_1+d_2+2}\times\dots\times\PP^{-d_{r-1}+d_r+2}$.

An algorithm can be elicited for determining which tuples $(U_1,\dots,U_n)$ can admit stable chains and which ones cannot.  Roughly, it works by recursion on rank.\benum\item Start with a tuple $(U_1,\dots,U_n)$ of rank $r-1$ and degree $d-a$, for some $a$, such that neither its slope nor the slopes of its subbundles exceeds or is equal to $d/r$.\item Replace $U_1$ of this chain with $U_1\plus\CO(a)$.\item Check if there is a subbundle of $E_a=(\CO(a)\plus U_1)\plus\cdots\plus U_n$ containing the $\CO(a)$ that has slope larger than or equal to $d/r$ and which is necessarily annihilated by all possible Higgs fields for $E_a$ that act with weight 1 on $(\CO(a)\plus U_1,\dots,U_n)$.\item If there is, discard $(\CO(a)\plus U_1,\dots,U_n)$.\item Repeat for the tuple $(\CO(a),U_1,\dots,U_n)$.\eenum

It can be shown that this algorithm terminates, as there are only finitely-many $a$ for which stability is possible, just as in the rank-2 case above.

For rank $3$ and degree $-1$, the ordered tuples of bundles admitting stable chains are $(\CO(1),\CO,\CO(-2))$, $(\CO(1),\CO(-1),\CO(-1))$, $(\CO,\CO,\CO(-1))$, and $(\CO\plus\CO,\CO(-1))$.  The first three are of type $(1,1,1)$ and the latter is of type $(2,1)$.  There are none of type $(3)$, as expected, but there are also none of type $(1,2)$ by the algorithm above.  The sets of chains on the $(1,1,1)$ tuples have Poincar\'e polynomials equal to those of $\PP^{-1+0+2}\times\PP^{-0-2+2}$, $\PP^{-1-1+2}\times\PP^{1-1+2}$, and $\PP^{-0+0+2}\times\PP^{0-1+2}$, respectively.  For $(\CO\plus\CO,\CO(-1))$, a map $\phi_1:\CO\plus\CO\rightarrow\CO(-1)\tensor\CO(2)$ is stable if and only if it is surjective.  If it is not surjective, then its image is either $0$, in which case the kernel is $\CO\plus\CO$ and therefore destabilizing, or is a sub-line bundle of degree $k<1$ in $\CO(-1)\tensor\CO(2)=\CO(1)$.  The kernel, accordingly, is a line bundle of degree $-k>-1$, which is destabilizing.  If the image is all of $\CO(1)$, then the kernel is isomorphic to $\CO(-1)$, and the resulting chain is stable.  Assuming now that $\phi_1$ is surjective, we have that the induced map $\widetilde\phi_1$ from global sections of $\CO\plus\CO$ to global sections of $\CO(1)$ must have full rank; that is, it must be an element of $\tb{GL}_2(\CC)$.  Quotienting by the right multiplication action of $\mbox{Aut}(\CO\plus\CO)=\tb{GL}_2(\CC)$ leaves only the identity, and so the set of chains on $(\CO\plus\CO,\CO(-1))$ has $b_0=1$ as its only nonzero Betti number.

What remains to be determined is the Morse index for each of these components of the fixed point set.  According to the formula \eqref{TheFormula}, the Morse index for $(\CO(1),\CO,\CO(-2))$ is $6$; for $(\CO(1),\CO(-1),\CO(-1))$, it is $4$; for $(\CO,\CO,\CO(-1))$, it is $2$; and for $(\CO\plus\CO,\CO(-1))$, it is $0$.  Putting all of this together, we get
\beqn\CP(3,-1;x) & = &  x^0(1)+x^2(1+x^2+x^4)(1+x^2)\notag
\\& &\quad\;+x^4(1+x^2+x^4)+x^6(1+x^2)\nonumber\eeqn\noin
which simplifies to $1+x^2+3x^4+4x^6+3x^8$.

It follows that the moduli space is topologically connected with three algebraic components.  That the Poincar\'e series is not palindromic indicates that, while the total space of $\CM(3,-1)$ is smooth, the nilpotent cone itself is not.\\

\begin{remark} There is no need to calculate the Betti numbers for $\CM(3,-2)$ because there is a degree duality taking points in $\CM(3,-1)$ to points in $\CM(3,-2)$, first by taking the dual co-Higgs bundle $(E^*,\phi^*)$, and then by tensoring $E^*$ by $\CO(-1)$.  On chains, this duality reverses the type of the chain, e.g. a $(1,1,1)$ chain goes to a $(1,1,1)$ chain, but a $(2,1)$ chain goes to a $(1,2)$ chain, and vice-versa.  This duality preserves the underlying topological structure of the moduli spaces.\\
\end{remark}

According to the algorithm, for rank $4$ and degree $-1$ the tuples admitting stable chains are
those given in Table~\hyperlink{tab:one}{\ref*{tabone}}.
There, we read a list of the form ``$\lst{a\,b \Div c \Div d\,e}$'' to mean \[(\CO(a)\plus\CO(b),\CO(c),\CO(d)\plus\CO(e)).\]  Computing Poincar\'e polynomials of sets of chains for the tuples in the table, and then combining the data with the Morse indices as in the rank-3 case, gives us
\beqn\CP(4,-1;x) &=& 1+x^2+3x^4+5x^6+9x^8+13x^{10}+18x^{12}+22x^{14}
\notag\\& &\quad\;+20x^{16}+10x^{18}.\nonumber\eeqn

\begin{remark} There are several tuples containing a rank-2 bundle, but there is no tuple with more than one rank-2 bundle.  Before the recent work of Garc\'ia-Prada, Heinloth, and Schmitt \cite{GHS:11}, the most formidable obstacle to computing the Betti numbers for the moduli space of ordinary rank-4 Higgs bundles was the existence of stable $(2,2)$ chains, which could not be directly attacked by Thaddeus' treatment of chains of length 2 \cite{MT:94}.\\
\end{remark}

\begin{remark} As with the rank-3 moduli space, there is no need to make a separate calculation for $\CM(4,-3)$, because of degree duality.\\
\end{remark}

\begin{table}[t] \hypertarget{tab:one}{\leavevmode}
\caption{Tuples admitting\label{tabone} stable chains for rank $4$ and degree $-1$.}
\noin\bcent\begin{math}
\begin{array}{|c|c|}
\hline
\mbox{\tb{Type}} & \mbox{\tb{Morse index}},[\mbox{\tb{Chain}}] \\
\hline\hline
(1,1,1,1) & 8,\,\lst{0\Div0\Div0\Div-1};\;\;8,\,\lst{0\Div1\Div-1\Div-1};\\& 10,\,\lst{0\Div1\Div0\Div-2};\;\;10,\,\lst{1\Div -1\Div 0\Div -1};\\ & 10,\,\lst{1\Div 0\Div -1\Div -1};\;\;12,\,\lst{1\Div 0\Div 0\Div -2};\\ & 12,\,\lst{1\Div 1\Div -1\Div -2};\;\;12,\,\lst{2\Div 0\Div -2\Div -1};\\ & 14,\,\lst{2\Div 0\Div -1\Div -2};\;\;16,\,\lst{2\Div 1\Div -1\Div -3}\\
\hline
(4),(3,1),(1,3),(2,2) & \mbox{no output}\\
\hline
(2,1,1) & 4,\,\lst{0\; 0\Div 0\Div -1};\;\;8,\,\lst{1\; 0\Div 0\Div -2}\\
\hline
(1,2,1) & 0,\,\lst{0\Div 0\; 0\Div -1};\;\;4,\,\lst{1\Div 0\; -1\Div -1}\\
\hline
(1,1,2) & 8,\,\lst{1\Div 0\Div -1\;-1}\\
\hline
\end{array}
\end{math}
\ecent
\end{table}

Finally, we note that the Poincar\'e polynomial for rank $5$ and degree $-1$ is
\begin{multline*} 1+x^2+3x^4+5x^6+10x^8+15x^{10}+26x^{12}+38x^{14}+56x^{16}+77x^{18}
\\+105x^{20}+131x^{22}+156x^{24}+165x^{26}+154x^{28}+103x^{30}+40x^{32}.\end{multline*}  
The calculations required for this result are markedly more difficult.  There are many more possible chains to contend with, there exist stable chains containing more than one rank-2 bundle, and there is a type-change phenomenon.  For the previous ranks, the holomorphic type of the bundles did not change within a component of the fixed point locus.  At rank 5, the holomorphic type of a $U_i$ may change within a component.

We also point out that while degree $-1$ and degree $-4$ necessarily have the same Betti numbers, degree $-2$ and $-3$ are a separate degree pair, unrelated to $-1$ or $-4$ by the duality mentioned earlier.

At rank $6$, there are several hundred admissible tuples of bundles.

\section{ADHM recursion formula}

In \cite{CDP:11}, Chuang, Diaconescu, and Pan give a recursion formula conjectured to relate the Donaldson--Thomas invariants of the usual Higgs bundle moduli space for genus $g\geq1$ to so-called ``asymptotic ADHM'' invariants.  In \cite{SM:12}, Mozgovoy finds a multivariable power series solution, and shows that the coefficients agree with the Hausel--Rodriguez-Villegas conjectures for Hodge polynomials of ordinary Higgs bundle moduli spaces \cite{HRV:08}.  Moreover, Mozgovoy solves a ``twisted'' version of the recursion formula and extends the solutions to genus 0.  These solutions can be conjectured to be Hodge polynomials of twisted Higgs bundles moduli spaces, where the Higgs field takes values in $\CO(t)$. In particular, for $g=0$ and $t=2$, these are the co-Higgs bundle moduli spaces.

For ranks 2 through 5, the conjectural Poincar\'e polynomials in \cite{SM:12} coincide with those in the previous section, therefore verifying conjectures presented in \cite{SM:12}.

Finally, we conjecture that the Betti numbers of co-Higgs moduli spaces on $\PP^1$ are independent of the degree.  This is implicit in the data coming from the ADHM formula: once the rank is fixed, there are no further parameters in the conjectural Poincar\'e polynomials.

Degree independence is known for ordinary Higgs bundles, but the proof uses properties of the character variety that are unavailable for co-Higgs bundles on $\PP^1$.  There is a diffeomorphism between the character variety of a higher genus curve and the moduli space of ordinary Higgs bundles on that curve, furnished by the nonabelian Hodge theorem originating in \cite{NJH:86,SKD:87,TSCH:92,KC:88,TSCH:88}, and the Riemann--Hilbert correspondence.  Degree independence of Betti numbers is proven for the character variety in \cite{HRV:08}.  Unfortunately, the nonabelian Hodge theorem depends in a crucial way on Higgs fields taking values in the canonical line bundle, and therefore does not extend in an obvious way to co-Higgs bundle moduli spaces on $\PP^1$.


\bibliographystyle{acm}

\end{document}